\documentclass{article}
\usepackage{amssymb}

\newtheorem{lemma}{Lemma}

\newtheorem{theorem}{Theorem}

\newtheorem{proposition}{Proposition}

\title{On a conjecture of Graham on the $p$-divisibility
of central binomial coefficients}

\author{Ernie Croot, Hamed Mousavi, Maxie Schmidt}

\begin{document}

\maketitle

\begin{abstract} We show that for every $r \geq 1$, and all
$r$ distinct (sufficiently large) primes 
$p_1,..., p_r > p_0(r,\varepsilon)$, there exist infinitely many integers $n$
such that ${2n \choose n}$ is divisible by these primes to only low
multiplicity.  From a theorem of Kummer, an upper bound for the 
number of times that a prime $p_j$ can divide ${2n \choose n}$ is
$1+\log n / \log p_j$; and our theorem shows that we can find 
integers $n$ where for $j=1,...,r$, $p_j$ divides ${2n \choose n}$
with multiplicity at most $\varepsilon$ times this amount.  
We connect this result
to a famous conjecture by R. L. Graham on whether there
are infinitely many integers $n$ such that ${2n \choose n}$
is coprime to $105$.
\end{abstract}

\section{Introduction}

In \cite{berend} and \cite{graham} it is mentioned that R. L. 
Graham had offered \$1,000 to settle the problem of whether
or not there are infinitely many integers $n$ such that
${2n \choose n}$ is relatively prime to $105 = 3 \cdot 5 \cdot 7$.
From the following theorem of Kummer \cite{kummer}
\bigskip

\noindent {\bf Kummer's Theorem:}\  
For a prime $p$ we have that the number of times that
$p$ divides ${n \choose m}$
equals the number of carries when adding the numbers $m$ and
$n-m$ in base-$p$.
\bigskip

\noindent we immediately see that Graham's problem is equivalent to asking
whether there are infinitely many integers $n\geq 1$ 
with the property
that when we add $n$ to itself in bases $3$, $5$, and $7$, there
are no carries.  In other words, are there infinitely many 
integers $n\geq 1$ such that all the base-$3$ digits are in
$\{0,1\}$, all the base-$5$ digits are in $\{0,1,2\}$, 
and all the base-$7$ digits are in $\{0,1,2,3\}$?  
If so, then there are 
infinitely many integers $n$ such that ${\rm gcd}\left ({2n \choose n}, 105\right )=1$; and if not, then there are at most finitely many 
integers $n\geq 1$ with ${\rm gcd}\left ({2n \choose n},105
\right )=1$.

In \cite{erdos}, Erd\H os, Graham, Ruzsa, and Straus proved that
for every {\it pair} of primes $p,q$, there are infinitely many 
integers $n\geq 1$ with ${\rm gcd}\left ( {2n \choose n}, pq\right )=1$; however, there are no such results in the literature for
$3$ or more primes (though, for example, there are results 
\cite{sanna, ford} on when ${2n \choose n}$ is coprime to $n$
and \cite{granville} when ${2n \choose n}$ is squarefree).  
Apart from whether one can give a proof of whether there are or aren't infinitely 
many $n$ with ${\rm gcd}\left ( {2n \choose n}, 105\right )=1$, 
one can at least ask whether it's
{\it plausible} or not that such integers $n \geq 1$ exist.  
Pomerance gave a simple heuristic for why there should exist
infinitely many $n\geq 1$ with this property (see, for example,
\cite{pomerance}):  if we choose a random $n \in [1,x]$, the 
probability that all its base-$3$ digits are in $\{0,1\}$ should be
about $(2/3)^{\log(x)/\log 3} \approx x^{-0.37}$; the probability that all its
base-$5$ digits are $\{0,1,2\}$ should be about 
$(3/5)^{\log(x)/\log 5} \approx x^{-0.32}$; and the probability that all its base-$7$ digits are $\{0,1,2,3\}$ 
should be about $(4/7)^{\log(x)/\log 7}\approx x^{-0.29}$.  
Assuming independence,
the probability that a random $n \in [1,x]$ satisfies all three
conditions is about $x^{-0.37}x^{-0.32}x^{-0.29} = x^{-0.98}$.  
So, we would expect there to be about $x^{0.02}$ numbers $n\in [1,x]$
with the property that ${\rm gcd}\left ( {2n \choose n}, 105\right )=1$, which clearly tends to infinity the larger we take $x$ to be.

One can extend Pomerance's heuristic to any number of odd primes,
making the same independence assumptions (that the events 
$E_1,...,E_r$ are mutually independent, 
where for a randomly chosen integer $n \in [1,x]$, $E_j$ is 
the event that the base-$p_j$ digits of 
$n$ are in $\{0,1,...,(p_j-1)/2\}$).  When one does this,
one would expect there to exist infinitely many integers $n \geq 1$
such that ${\rm gcd}\left ( {2n \choose n}, p_1 \cdots p_r\right ) = 1$, for distinct odd primes $p_1,...,p_r$, provided that
\begin{equation}\label{condition}
- \sum_{j=1}^r {\log \left ( {1 \over 2} + {1 \over 2 p_j} \right )
\over  \log(p_j)}\ <\ 1;
\end{equation}
and that (using the Borel-Cantelli Lemma) 
there should be only {\it finitely} many such $n$ if the $>$ is
replaced with a $<$.  We make no guesses about the possible case
when the left-hand-side equals $1$, exactly -- if it is even 
possible.  See \cite[example 2.1]{berend2} for related discussion.

What is interesting here is that even if we consider a slight 
weakening of the problem where we allow ${2n \choose n}$ to 
be divisible by the primes $p_1,...,p_r$ to {\it low multiplicity}, we
get the same condition (\ref{condition}) guaranteeing the existence
of infinitely many such $n\geq 1$:  in light of Kummer's
theorem, the number of times that a prime $p_j$ can divide a number
$n$ is at most about $1+\log(n)/\log p_j$, 
since this is an upper bound on the 
number of base-$p_j$ digits of $n$.    
If we select a random 
$n \in [1,x]$, the probability that all but at most 
$k$ of the base-$p_j$ digits are in $\{0,1,2,...,(p_j-1)/2\}$ 
is 
$$
\asymp\ {[\log x / \log p_j] \choose k} \left ( {1 \over 2} + {1 \over 2 p_j} \right )^{\log(x)/\log p_j - k} \left ( {1 \over 2} - {1\over 2 p_j}
\right )^k,
$$
for $k = o(\log x)$.  This has size (assuming $k=o(\log x)$)
$$
\left ( {1 \over 2} + {1 \over 2 p_j}\right )^{(1-o(1))\log(x)/\log p_j},
$$
which, apart from the factor $1-o(1)$ in the exponent, has the
same form as the probability for the case where {\it every} 
base-$p_j$ digit of $n$ is in $\{0,1,...,(p_j-1)/2\}$.  Making
the same independence assumptions as before, we thus would
expect that if (\ref{condition}) holds, then there should 
exist infinitely many integers $n\geq 1$ where for 
$j=1,...,r$, $p_j$ divides ${2n \choose n}$ to multiplicity at
most $o(\log n)$; and, if, instead, the left-hand-side of 
(\ref{condition}) is $> 1$, we would expect there to be only
{\it finitely} many such $n\geq 1$. 
\bigskip

In this paper, we don't quite prove that (\ref{condition}) implies
there are infinitely many such $n\geq 1$, but we do prove something
in this direction:

\begin{theorem}\label{maintheorem}  Suppose $r \geq 1$, 
$0 < \varepsilon < 1/20r^2$, and let $p_1,...,p_r \geq c_0(r,\varepsilon)$ be distinct primes, where $c_0(r,\varepsilon)$ is
some function of $r$ and $\varepsilon$ (can be deduced from the proof).  Then, there is a sequence
$n_1, n_2, ...$ of integers $n$ such that for all $i=1,...,r$,
$$
\nu_{p_i}\left ( {2n \choose n} \right )\ \leq\ 
{\varepsilon \log n \over \log p_i},
$$
where $\nu_p(x)$ denotes the number of times the prime $p$ divides 
$x$.
\end{theorem}

As we said, a trivial upper bound for 
$\nu_{p_i}( {2n \choose n})$ is $1 + (\log n)/\log p_i$, since
$n$ has at most this many base-$p_i$ digits; so the theorem is
saying that we can find infinitely many $n$ where we are smaller
than this amount by a factor $\varepsilon$, for all the primes
$p_1,...,p_r$, simultaneously.

As one will see, the proof is fairly technical.  What would greatly
simplify it is if one had that the numbers 
$1/\log 2, 1/\log p_1, ..., 1/\log p_r$ were linearly 
independent over the rationals.  This is not known to be true
for arbitrary sets of primes, but it would follow from the
following conjecture:
\bigskip

\noindent {\bf Schanuel's Conjecture \cite{lang}}\  Given any $n$ complex
numbers $z_1,...,z_n$ that are linearly independent over the
rationals, the field extension ${\mathbb Q}(z_1,...,z_n,e^{z_1},...,
e^{z_n})$ has transcendence degree at least $n$ over ${\mathbb Q}$.
\bigskip

If Schanuel's Conjecture holds, then taking $n=r+1$,
and taking $z_1=\log 2$, $z_2 = \log p_1$, ..., $z_{r+1} = \log p_r$,
we see that ${\mathbb Q}(\log 2, \log p_1, ..., \log p_r)$ has
transcendence degree $r+1$.  Now suppose we had a linear combination
$$
{\lambda_1 \over \log 2} + {\lambda_2 \over \log p_1} + \cdots + 
{\lambda_{r+1}\over \log p_r}\ =\ 0,
$$
where\ $\lambda_1,...,\lambda_{r+1} \in {\mathbb Q}$ and not all $0$.
Without loss, assume that $\lambda_1 \neq 0$.  Then, the linear
relation would imply that 
$$
{\mathbb Q}(\log 2, \log p_1, ..., \log p_r)\ =\ 
{\mathbb Q}(\log p_1, \log p_2, ...,\log p_r),
$$
which can have transcendence degree at most $r$, which would be
a contradiction.  Thus, no such linear relations can hold.

\section{Proof of the Main Theorem}

For a real number $x \geq 0$ we will use the notation $\{x\} = x - [x]$ to denote the fractional part of $x$, and where $[x]$ denotes the integer part of $x$.  We will also use $\| x\|$ to denote the distance from $x$ to the nearest integer.
\bigskip

As we said in the introduction, the 
central binomial coefficients in the statement of
the theorem are somewhat of a distraction, in light of Kummer's
Theorem.  This theorem implies that for $p_j$ odd, if all but at most
$\varepsilon (\log n)/\log p_j$ of the base-$p_j$ digits of $n$ are $\leq (p_j-1)/2$,
then $v_{p_j}\left ( {2n \choose n}\right ) \leq 
\varepsilon (\log n)/\log p_j$; and
establishing that there are infinitely many integers $n$ with
this property (few base-$p_j$ digits that are $> (p_j-1)/2$) 
is the path we will take to prove Theorem \ref{maintheorem}.
\bigskip

In carrying out this verification, we will 
make use of the following theorem:

\begin{theorem}\label{theorem2}  Suppose that $p_1,...,p_r$ are distinct odd primes.  For $i=1,...,r$, and $n \geq 1$, define
$$
\alpha_i(n)\ :=\ p_i^{\{n (\log 2)/ \log p_i\}-1}\ =\ 
p_i^{n(\log 2)/\log p_i - [n(\log 2)/\log p_i] - 1},
$$
Define, for $H \geq 1$ and $i=1,2,...,r$,
\begin{equation}\label{Uidef}
U_i(H)\ :=\ \left \{ {d_1 \over p_i} + {d_2 \over p_i^2} + 
\cdots + {d_H \over p_i^H}\ :\ 0\ \leq\ d_1, ..., d_H\ \leq\ 
{p_i \over 3}\right \}\ +\ \left [ 0,\ {1 \over p_i^H} \right ).
\end{equation}
Then, we have
that for $N \geq 1$, for $H$ tending slowly to infinity with $N$
in a way that can be made precise by following the proof, and for arbitrary sequences of real numbers
$\{\beta_i(n)\}_{n=1}^\infty$, $i=1,...,r$, 
\begin{equation}\label{conclusion0}
{\#\{ n \leq N\ :\ \exists\ s \leq 10^{10r^2H}\ \forall\ j =1,...,r,\ 
\{ s \alpha_j(n) + \beta_j(n)\} \in U_j(H)\} \over N}\ \geq\ 1-o(1).
\end{equation}
\end{theorem}

\noindent {\bf Remark:}  We could prove a theorem like this where $H$
is {\it bounded}, but the price to be paid is that the error term $o(1)$
on the right-hand-side of (\ref{conclusion0}) would have to be some function of $H$ that tends to $0$ as $H \to \infty$.  It's simpler just to take $H$ tending to infinity than bother with expressing the error term in terms of $H$.
\bigskip

Now let us see that this theorem implies Theorem \ref{maintheorem}:  it clearly suffices to prove that for each integer $N$ sufficiently large, we can find an integer $n$ satisfying
$$
2^{N/2}\ <\ n\ \leq\ 2^N,
$$
so that for all $j=1,2,...,r$, all but at most 
$\varepsilon (\log n)/\log p_j$ of the 
base-$p_j$ digits of $n$ are $\leq p_j/3$.  

So, let us suppose $N$ is given, and suppose $H = H(N)$ be some
function slowly tending to infinity with $N$ as in the above theorem
($H(N) \asymp \log\log N$ would be sufficient for our purposes).  Let $p$ be the smallest of the primes
$p_1,...,p_r$, and then define $\ell$ to be the unique integer such that
$$
4^\ell\ <\ p^H\ <\ 4^{\ell+1}.
$$
Then, from Theorem \ref{theorem2} and the pigeonhole principle, for some $t=0,1,2,...,\ell-1$ we will have that for all but $o(N/\ell)$ of the 
integers $0 \leq m < N/\ell$ the following holds:  for any
$\beta \in {\mathbb R}^r$, there exists $1 \leq s \leq 10^{10r^2 H}$
such that for every $j=1,...,r$,
$$
\left \{ s \alpha_j(m\ell + t) + \beta_j \right \}
\ \in\ U_j(H).
$$

Let 
$$
N'\ :=\ \lfloor N/\ell\rfloor.
$$
We now construct the number 
$$
n\ :=\ n_0 2^{\ell N'+t} + n_1 2^{\ell (N'-1)+t} + \cdots + n_{N'}2^t
$$
as follows:  we start by letting $n_0=1$.  Assume we have constructed $n_0,..., n_{d-1}$.  Now we show how to 
construct $n_d$, $d \geq 1$:  for $j=1,...,r$, we let
$$
\beta_j(\ell (N'-d)+t)\ :=\ 
{n_0 2^{\ell N'+t} + n_1 2^{\ell (N'-1)+t} +\cdots + n_{d-1} 2^{\ell (N' - d+1)+t}\over p_j^{m_{j,d}}},
$$
where for an integer $h$ we define
$$
m_{j,h}\ =\ \left [ {(\ell (N'-h) + t) \log 2 \over \log p_j}\right ] + 1.
$$
(Alternatively:  $m_{j,h}$ is the unique integer so that 
$2^{\ell(N'-h)+t} / p_j^{m_{j,h}}$ lies in $[1/p_j, 1)$.)  

If it exists, we let $1 \leq n_d \leq 10^{10r^2 H}$ be any integer where
\begin{equation}\label{Uhold}
\{ n_d \alpha_j(\ell(N'-d)+t) + \beta_j(\ell(N'-d)+t)\}\ \in 
U_j(H).
\end{equation}
If no such $n_d$ exists, just let $n_d=0$.
\bigskip

In order to see that this construction works, we begin by noting that
for any integer $h$,
$$
\alpha_j(h)\ =\ p_j^{\{h (\log 2)/\log p_j\}-1}\ =\ 
{p_j^{h (\log 2)/\log p_j}\over p_j^{[h (\log 2)/\log p_j]+1}}\ =\ 
{2^h \over p_j^{h'}},
$$
where $h'$ is the unique integer such that this belongs to the interval
$[1/p_j, 1)$.  

Thus, when we go to construct $n_d$, we will have
$$
n_d\alpha_j(\ell(N'-d)+t) + \beta_j(\ell(N'-d)+t)
\ =\ {n_0 2^{\ell N'+t} + n_1 2^{\ell (N'-1)+t} +\cdots + n_d 2^{\ell(N'-d)+t} \over p_j^{m_{j,d}}}.
$$

It follows that if we write
\begin{equation}\label{2pjsum}
n_02^{\ell N'+t} + n_1 2^{\ell(N'-1)+t} + \cdots + n_d 
2^{\ell(N'-d)+t}\ =\ c_0 + c_1 p_j + c_2 p_j^2 + \cdots,
\end{equation}
where $0 \leq c_i \leq p_j-1$, then from (\ref{Uhold}) we 
deduce that if $n_d \neq 0$ then 
$$
0\ \leq\ c_{m_{j,d}-1},\ c_{m_{j,d}-2},\ ...,\ c_{m_{j,d}-H}
\ \leq\ {p_j \over 3},
$$
and so in particular, since $m_{j,d} - m_{j,d+1} < H$,
we have that $0 \leq c_u \leq p_j/3$ for 
$$
m_{j,d+1}\ \leq\ u\ \leq\ m_{j,d}.
$$

Now, if we continue adding on additional terms to (\ref{2pjsum}),
\begin{equation}\label{extraterms}
n_{d+1} 2^{\ell(N'-d-1)},\ n_{d+2} 2^{\ell(N'-d-2)},\ ...
\end{equation}
these will only have an effect on the terms $c_z p_j^z$ where 
$$
z\ \leq\ m_{j,d+1} + [(\log n_{d+1})/\log p_j]+1\ \leq\ 
m_{j,d+1} + 10 r^2 H (\log 10)/\log p_j + 1.
$$
Thus, the terms $c_u p_j^u$ where
$$
m_{j,d+1}+ 10 r^2 H(\log 10)/\log p_j + 1\ <\ u\ \leq\ m_{j,d}
$$
in (\ref{2pjsum}) will be unchanged, as will all the other higher-order terms with $u > m_{j,d}$.  
\bigskip

Now we distinguish two possibilities for each $d \leq N'$:  
we let $D^\sharp$ be those $d$ such that there {\it does} 
exist an $n_d \leq 10^{10 r^2 H}$ where (\ref{Uhold})
holds, and we let $D^\flat$ be those $d$ for which it doesn't.
Recall that Theorem \ref{theorem2} and the pigeonhole principle
were used to show that the number of such 
$d$ where (\ref{Uhold}) doesn't hold is 
$o(N/\ell) = o(N')$; and so, $|D^\flat| = o(N')$.  So,
$|D^\sharp| = N'(1-o(1))$.

For each $d \in D^\sharp$ we have that for each $j=1,...,r$, 
at most $10r^2 H (\log 10)/\log p_j + 1$ base-$p_j$ digits 
$c_u$ with $m_{j, d+1} \leq u \leq m_{j,d}$ are $> p_j/3$.
So, the $d\in D^\sharp$ contribute a total of at most
\begin{eqnarray*}
|D^\sharp| (10 r^2 H (\log 10)/\log p_j + 1)\ &\ll&\ 
N' r^2 H/\log p_j \\
&\ll&\ N r^2 / (\log p)(\log p_j)
\end{eqnarray*}
base $p_j$ digits $> p_j/3$.
If the smallest of the primes $p_1,...,p_r$ is large enough then for all $j=1,...,r$ we will have that
this can be made smaller than 
$\varepsilon N/2\log p_j$, say.

And for each $d \in D^\flat$, in the worst cast for every $j=1,...,r$,
all of the $c_u$ with $m_{j,d+1} \leq u \leq m_{j,d}$ could be $> p_j/3$.  Note that in this case (the case $d \in D^\flat$)
there are at most $\ell+1$ bad digits $c_u$ with $m_{j,d+1} \leq u
\leq m_{j,d}$.

All told, then, the total number of bad 
base-$p_j$ digits that are $>p_j/3$ in this case, over all
$d \in D^\flat$, is at most
$$
(\ell + 1) o(N')\ =\ o(N),
$$
for every $j=1,...,r$.

In total, then, for every $j=1,...,r$, the number of bad 
base-$p_j$ digits (that are $> p_j/2$) is
at most 
$$
\varepsilon N / \log p_j.
$$
This is just what we need to show in order to prove
Theorem \ref{maintheorem}.

\section{Proof of Theorem \ref{theorem2}}

In proving this theorem we will need to understand how the 
vectors
\begin{eqnarray*}
&& (\alpha_1(n),\ \alpha_2(n),\ ...,\ \alpha_r(n))\\
&& \ \ \ \ \ 
=\ (p_1^{\{ n (\log 2)/\log p_1\} - 1},\ p_2^{\{n (\log 2) / \log p_2\}-1},\ ..., p_r^{\{n(\log 2)/\log p_r\}-1}).
\end{eqnarray*}
are distributed, as we vary over $n \leq N$.  

\subsection{The 2-dimensional case}

To better understand what is going on, we first consider the case 
where $r=2$.  There are two possibilities:  the first possibility
is that there do not exist integers $n_1, n_2, n_3$, with 
$n_1, n_2, n_3 \neq 0$, such that
\begin{equation}\label{linearcomb}
n_1 {\log 2 \over \log p_1} + n_2 {\log 2 \over \log p_2}\ =\ n_3.
\end{equation}

If this occurs, then as a consequence of 

\begin{theorem}[Multidimensional Weyl's Theorem]\label{weyl} Suppose $1,\vartheta_1,...,\vartheta_r$ are real numbers that are linearly
independent over the rationals.  Then, for 
$\vec \vartheta = (\vartheta_1, ..., \vartheta_r) \in {\mathbb R}^r$,
the sequence $\{k \vec \vartheta\}_{k=1}^\infty$ is 
uniformly distributed in ${\mathbb R}^r / {\mathbb Z}^r$.
\end{theorem}
(see \cite[example 6.1]{kupiers}) we have that the vector 
$(n (\log 2)/\log p_1, n (\log 2)/\log p_2)$ is uniformly
distributed mod $1$ as we vary over $n =1,2,3,...$; and, therefore,
the set $(\alpha_1(n), \alpha_2(n))$ is dense in the box
$[1/p_1, 1]\times [1/p_2, 1]$.  

The second possibility is that there {\it do} exist integers 
$n_1,n_2,n_3 \neq 0$ such that (\ref{linearcomb}) holds
(If $n_1$ were allowed to be $0$, then we would have that
(\ref{linearcomb}) implies $n_2 \log 2 = n_3 \log p_2$, which can
only hold if $n_2 = n_3 = 0$; and a similar thing occurs for when 
$n_2 = 0$ or when $n_3 = 0$; so, if one of these $n_i$ were $0$, the others would have to be as well.)

By multiplying through by $-1$ as needed, we can assume $n_2 > 0$; and we will
assume that the $p_1$ and $p_2$ are arranged so that 
$$
|n_1/n_2|\ \leq\ 1.
$$

We will show that the set 
\begin{equation}\label{vectors}
(\alpha_1(n),\ \alpha_2(n))\ =\ (p_1^{\{ n (\log 2)/\log p_1\}-1},
\ p_2^{\{n (\log 2)/\log p_2\}-1}),\ n=1,2,3,...
\end{equation}
is contained in a union of a finite set of non-linear curves.  It turns out that, moreover, the set is equidistributed on these curves (when
we restrict to $[1/p_1,1]\times [1/p_2,1]$) with respect to the 
right measure; though, we don't actually need the full strength of such a statement, so don't bother to prove it.
\bigskip

We claim that for each integer $n \geq 1$,
$$
\left \{ {n \log 2 \over \log p_2}\right \}\ =\ f(n) - {n_1 \over n_2} \left \{ {n \log 2 \over \log p_1}
\right \},
$$
where $f(n) \in S$, a finite set of possibilities.
\footnote{Note that we have redefined the function $f$ from how we used it in a previous section.}
To see this,
we begin by rewriting 
(\ref{linearcomb}) as
\begin{equation}\label{linearcomb2}
{n \log 2 \over \log p_2}\ =\ {n n_3 \over n_2} - {n_1\over n_2} {n \log 2 \over \log p_1}.
\end{equation}
Now we write
\begin{equation}\label{that1}
{n \log 2 \over \log p_1}\ =\ \left [ {n \log 2 \over \log p_1}\right ]\ +\ \left \{ {n \log 2 \over \log p_1}\right \}\ =\ k(n) n_2 + a(n) + \left \{ {n \log 2 \over \log p_1}\right \},
\end{equation}
where $k(n)$ is an integer, and $0 \leq a(n) \leq n_2-1$.  We also write
$$
n\ =\ \ell(n) n_2 + b(n),\ {\rm where\ }0\leq b(n) \leq n_2-1,\ {\rm and\ } \ell(n) \in {\mathbb Z}.
$$
It follows, then, upon plugging this and (\ref{that1}) into (\ref{linearcomb2}), that
$$
{n \log 2 \over \log p_2}\ =\ \ell(n) n_3+ {b(n) n_3 \over n_2} - k(n) n_1 - {a(n) n_1 \over n_2} - {n_1\over n_2}
\left \{ {n \log 2 \over \log p_1} \right \}.
$$
Thus, since $|n_1/n_2| \leq 1$,
$$
\left \{ {n \log 2 \over \log p_2}\right \}\ =\ \left \{ {b(n) n_3 \over n_2} - {a(n) n_1 \over n_2} \right \}
- {n_1 \over n_2} \left \{ {n \log 2 \over \log p_1} \right \} + \delta,\ {\rm where\ }\delta \in 
\{0,1,-1\}.
$$
We would take $\delta=0$ if the preceding 
terms add to a number in $[0,1)$; take $\delta=1$ if they produce a number in $[-1,0)$; and take $\delta=-1$
if they produce a number in $[1,2]$.

It follows that we may take $S$ to be
$$
S\ =\ \left \{ \left \{ {b n_3\over n_2} - {a n_1 \over n_2}\right \}\ :\ 
a,b=0,1,...,n_2-1\right \}\ +\ \{0,1,-1\}.
$$
Thus,
$$
|S|\ \leq\ 3 n_2^2.
$$

We conclude that, for $n \geq 1$,  
\begin{eqnarray}\label{points}
&&(p_1^{ \{ n \log 2 / \log p_1\}-1},\ 
p_2^{\{ n \log 2 / \log p_2\}-1})\nonumber\\ 
&&\hskip0.5in=\ 
(p_1^{ \{ n \log p / \log p_1\}-1},\ 
c(n) p_2^{-(n_1/n_2) \{n \log p / \log p_1\}-1}),
\end{eqnarray}
where $c(n) = p_2^{f(n)}$, where, recall, $f(n) \in S$.

As we vary over $n \leq N$, and let $N \to \infty$, 
all the points (\ref{points}) lie on a set of at most $|S| \leq 3n_2^2$ curves of the form
$$
C_s\ :=\ \{ (p_1^{t-1},\ c_s p_2^{ -(n_1/n_2) t-1}\ :\ 0 \leq t < 1\},\ {\rm where\ }c_s = p_2^{s},\ {\rm 
where\ }s\in S.
$$

\subsection{None of these curves are lines}

These curves are just dilates of one another in the second coordinate.  So, to show
that none are lines, it suffices to show that the curve with points
$$
z(t)\ :=\ (p_1^{t},\ p_2^{ -(n_1/n_2)t}),
$$
is not a line.

To see this it suffices to prove that 
$$
p_1\ \neq\ p_2^{-(n_1/n_2)},
$$
which is clearly the case, since upon raising both sides to the $n_2$ power,
if they were equal we would have
$$
p_1^{n_2}\ =\ p_2^{-n_1},
$$
which can't hold if $p_1$ and $p_2$ are distinct primes.

\subsection{Generalizing to higher dimensions} \label{generalsection}

Now suppose we have $r$ primes $p_1, ..., p_r$, and we wish to 
understand the possible vectors
\begin{equation}\label{pointsn}
(p_1^{\{ n \log 2/\log p_1\}-1},\ p_2^{\{n \log 2 /\log p_2\}-1},\ ...,\ 
p_r^{\{n \log 2 / \log p_r\}-1}),
\end{equation}
given that we have relations similar to (\ref{linearcomb}).  In this case, there can be 
more than one such relation.  We can express this set of relations as
\begin{eqnarray*}
a_{1,1} {\log 2 \over \log p_1} + a_{1,2} {\log 2 \over \log p_2} + \cdots + 
a_{1,r} {\log 2 \over \log p_r}\ &=&\ a_{1,r+1} \\
 &\vdots& \\
a_{k,1} {\log 2 \over \log p_1} + a_{k,2} {\log 2 \over \log p_2} + \cdots + 
a_{k,r} {\log 2 \over \log p_r}\ &=&\ a_{k,r+1},
\end{eqnarray*}
where all the $a_{i,j} \in {\mathbb Q}$, where $k \leq r-1$, and where
all these relations are linearly independent.  Note that
if there were $k=r$ linearly independent relations, then this would imply
that all the $\log 2 / \log p_i$ are rational numbers, which would imply
that for each $i=1,...,r$, $\log 2$ and $\log p_i$ are linearly dependent over the rationals, 
which we know is false, as it would imply that there is an integer
power of $2$ that equals an integer power of $p_i$.

Upon applying row-reduction to these equations, and permuting the $p_j$'s as needed, we
can reduce the above system to the following one:  for $j=1,...,k$, we have
$$
{\log 2 \over \log p_{r-j+1}}\ =\ 
b_{j,1} {\log 2 \over \log p_1} + \cdots + b_{j,r-k} {\log 2 \over \log p_{r-k}} + b_{j,r+1},
$$
where the $b_{j,h} \in {\mathbb Q}$.  We have, also, that (recalling
that the $p_j$'s have been permuted from their original ordering)
\begin{equation}\label{independentlogs}
1,\ {\log 2 \over \log p_1},\ ...,\ {\log 2 \over \log p_{r-k}}
\ {\rm are\ independent\ over\ }{\mathbb Q}.
\end{equation}
We note that this holds also in the case $k=0$, 
where there are {\it no} linear relations as above.

Getting a common denominator, we can rewrite the above as:  
for $j=1,...,k$, we have
\begin{equation}\label{nlogp0}
{\log 2 \over \log p_{r-j+1}}\ =\ 
{m_{j,1}\over n_j} {\log 2 \over \log p_1} + \cdots + {m_{j,r-k}\over n_j} 
{\log 2 \over \log p_{r-k}} + {m_{j,r+1}\over n_j},
\end{equation}
where, for all $j=1,...,k$ and $h=1,...,r-k,r+1$, the $n_j \geq 1$ and the $m_{j,h}$ are integers.

Now, we claim that for $n \geq 1$, and $j=1,...,k$,
\begin{equation}\label{nlogp}
\left \{ {n \log 2 \over \log p_{r-j+1}}\right \}\ =\ g_j(n)\ +\ {m_{j,1}\over n_j}
\left \{ {n \log 2 \over \log p_1}\right \} + \cdots + {m_{j,r-k} \over n_j}
\left \{ {n \log 2 \over \log p_{r-k}}\right \},
\end{equation}
where $g_j(n)$ takes on values in a finite set $S$ of possibilities.  

To see this, we proceed as with the 2-dimensional case:  for $j=1,...,k$ and
$h=1,...,r-k$, we define the numbers
$\ell_{j,h}(n) \in {\mathbb Z}$ and $0 \leq a_{j,h}(n) \leq n_j-1$ as follows
$$
{n \log 2 \over \log p_h}\ =\ \ell_{j,h}(n) \cdot 
n_j + a_{j,h}(n)\ +\ \left \{ {n \log 2 \over \log p_h} \right \}.
$$
Thus, from (\ref{nlogp0}) we have that
$$
{n \log 2 \over \log p_{r-j+1}}\ =\ \sum_{h=1}^{r-k} \ell_{j,h}(n) m_{j,h} + {a_{j,h}(n)m_{j,h}
\over n_j} + {m_{j,h} \over n_j} \left \{ {n \log 2 \over \log p_h} \right \} + n b_{j,r+1}.
$$
Taking the fractional part of both sides, we find that
\begin{equation}\label{values}
\left \{ {n \log 2 \over \log p_{r-j+1}}\right \}\ =\ 
\left \{ \sum_{h=1}^{r-k} {a_{j,h}(n) m_{j,h} \over n_j} \right \} + 
\left ( \sum_{h=1}^{r-k} {m_{j,h} \over n_j} 
\left \{ {n \log 2 \over \log p_h}\right \}\right ) + 
\{n b_{j,r+1}\} + \delta_j,
\end{equation}
where $\delta_j$ is an integer chosen so as to make the right-hand-side of this equation
be a real number in $[0,1)$.  Clearly, 
$\delta_j\ \in\ \{-\Delta, -\Delta+1, ..., 0, ..., \Delta\}$, where
$$
\Delta\ =\ 1 + r\cdot \max_{j,h} \left \lfloor 
{|m_{j,h}| \over |n_j|} \right \rfloor.
$$
Thus, if we let
$$
S\ :=\ \{c_j / n_j\ :\ j=1,...,k,\ {\rm and\ }0 \leq c_j \leq n_j-1\}\ +\ \{-\Delta, -\Delta+1, ..., 0,\ 1, ...,\ \Delta\},
$$
then from (\ref{values}) we see that
$$
\left \{ {n \log 2 \over \log p_{r-j+1}}\right \}\ =\ g_j(n)\ + \sum_{h=1}^{r-k} {m_{j,h} \over n_j} 
\left \{ {n \log 2 \over \log p_h}\right \},
$$
where $g_j(n) \in S$.

Thus, proceeding as in the 2-dimensional case, we see that the set of points (\ref{pointsn}) all
lie on one of the following finite set of surfaces given as follows:
\begin{equation}\label{surfaceform}
(p_1^{t_1-1},\ p_2^{t_2-1},\ ...,\ p_{r-k}^{t_{r-k}-1},\ c_1 p_{r-k+1}^{\theta_1(t_1,...,t_{r-k})-1},\ 
...,\ c_k p_r^{\theta_k(t_1,...,t_{r-k})-1}),
\end{equation}
where 
$$
{\rm for\ }i=1,...,k,\ 
c_i\ =\ p_{r-k+i}^{s_i},\ {\rm for\ some\ }s_i \in S,
$$
and where
$$
\theta_i(t_1, ..., t_{r-k})\ =\ \sum_{h=1}^{r-k} {m_{k-i+1, h} \over n_{k-i+1}} t_h.
$$
We note that for $k=0$ (no linear relations) the surface
(\ref{surfaceform}) just becomes 
$$
(p_1^{t_1-1},\ p_2^{t_2-1},\ ...,\ p_r^{t_r-1}).
$$

\subsection{Passing to parameterized curves} \label{paramsection}

\subsubsection{An illustrative example}

We would like to break these surfaces up into a union of parameterized curves of the form
$$
(c_1\cdot \alpha_1^{t}, c_2\cdot \alpha_2^{t}, ..., c_r\cdot \alpha_r^{t}),
$$
where the $\alpha_j$'s are all distinct, and none of the $c_j$'s are $0$.  One reason for doing this is that it will make certain estimates later on notationally simpler; another reason is that working with general surfaces would require proving more general versions of certain intermedial results (like Lemma 1); a third reason is that we don't rule out the 1-dimensional curve case actually happening, so a general approach that uses arbitrary surfaces would have to plan for it and include it as a special case, anyways.  We could just keep everything as one large surface, in the general case, and then later when we apply Fourier analysis, pay the price of integrating over (or summing over) that surface -- of variable dimension -- and then worry about points where certain Fourier transforms are large in magnitude.  However, once we have the surface as a union of curves, that stage of the argument only needs a simple sum over a single 1-dimensional variable.

One attempt at doing this
would be to take a surface of the form (\ref{surfaceform}), and set all but one of the $t_j$'s
to fixed values.  For example, if our surface were of the form
\begin{equation}\label{testsurface}
(p_1^{t_1},\ p_2^{t_2},\ p_3^{t_1+t_2}),\ t_1, t_2 \in [0,1),
\end{equation}
then if we were to freeze $t_1$ and vary $t_2$, we would get a curve of the form
$$
(c,\ p_2^{t},\ d\cdot p_3^{t}).
$$
Unfortunately, only two coordinates vary, not all three.  However, if we parameterize 
differently, then we can get a full, 3-dimensional curve:  let $t_1 = t$, 
$t_2 = t+\delta$ mod $1$, where $\delta$ is fixed and $t$ varies in $[0,1)$.  Then,
we get the parameterized curve
\begin{equation}\label{examplecurve}
(p_1^{t},\ c\cdot p_2^{t},\ d\cdot p_3^{2t}),
\end{equation}
where $c = p_2^{\delta}$, $d = p_3^{\delta}$.  Actually, this isn't quite right, since, for
example $2t > 1$ for $t > 1/2$; we need to introduce another curve to account for these 
possibilities.  Basically, we consider all curves of (\ref{examplecurve}) where
$c \in \{p_2^{\delta}, p_2^{\delta-1}\}$ and $d \in \{p_3^{\delta}, p_3^{\delta-1},
p_3^{\delta-2}\}$.  This covers all the cases; and, as we vary over all $\delta \in [0,1)$, we get a union of curves, where this union is exactly
the surface (\ref{testsurface}) when restricted to $(1/p_1,1] \times (1/p_2,1] \times (1/p_3,1]$.  

\subsubsection{Applying this idea to the surface (\ref{surfaceform})} \label{paramcurvesection}

To attempt something similar for the surfaces (\ref{surfaceform}), we will choose
\begin{equation}\label{ts}
t_1 = t,\ t_2 = L t + \delta_1,\ t_3 = L^2t + \delta_2,\ ...,\ t_{r-k} = L^{r-k-1}t + \delta_{r-k-1},
\end{equation}
where the $\delta_1,...,\delta_{r-k-1} \in [0,1)$, and where $L$ is an integer chosen suitably large
so that 
$$
\rho_i\ :=\ \theta_i(1,L,L^2,...,L^{r-k-1})\ \neq\ 0.
$$
It isn't hard to see that one can take 
$$
L\ \leq\ 2\cdot {\rm lcm}(n_1, ..., n_k)\max_{i,j} |m_{i,j}|.
$$

Using the parameterization (\ref{ts}), we will get that
$$
\theta_i(t_1,...,t_{r-k})\ =\ t\cdot \theta_i(1,L,L^2, ..., L^{r-k-1}) + \theta_i(0,\delta_1, ..., \delta_{r-k-1})\ =\ t\cdot \rho_i + \theta_i(0,\delta_1,...,\delta_{r-k-1})
$$
and applying this to (\ref{surfaceform}), we will get curves of the form
\begin{equation}\label{newcurves}
(p_1^{t-1}, d_2 p_2^{Lt-1},\ ...,\ d_{r-k} p_{r-k}^{L^{r-k-1}t-1},\ d_{r-k+1} p_{r-k+1}^{\rho_1 t-1},\ 
...\ d_r p_r^{\rho_k t-1}),
\end{equation}
where 
$$
d_2 = p_2^{\delta_1},\ d_3 = p_3^{\delta_2},\ ...,\ d_{r-k} = p_{r-k}^{\delta_{r-k-1}},
$$
and
$$
d_{r-k+1} = p_{r-k+1}^{\theta_1(0,\delta_1,...,\delta_{r-k-1})} c_1,\ ...,\ 
d_r = p_r^{\theta_k(0,\delta_1, ..., \delta_{r-k-1})} c_k,
$$
where the $c_i$ are of the form $p_{r-k+i}^{s_i}$, where $s_i \in S$.

Similar to how we dealt with (\ref{examplecurve}) not quite covering all possible curves, we actually need to
expand the set of possibilities for the $d_i$, given a fixed choice for $\delta_1,...,\delta_{r-k}$
(that is, our curves (\ref{newcurves}) don't quite cover everything):  we need to
also include dilates by integral powers of the $p_i$, to handle, for example,
$p_2^{Lt}$ not always being in the range $(1/p_2,1]$ (basically, the exponent $Lt$ 
needs to be considered mod $1$).  Thus, in fact, we want to consider the $d_i$'s in
dilated sets
\begin{equation}\label{dchoice1}
d_2 \in p_2^{\delta_1} D_2,\ d_3 \in P_3^{\delta_2} D_3,\ ...,\ d_{r-k} \in 
p_{r-k}^{\delta_{r-k-1}} D_{r-k},
\end{equation}
and
\begin{equation}\label{dchoice2}
d_{r-k+1} \in p_{r-k+1}^{\theta_1(0,\delta_1,...,\delta_{r-k-1})} c_1D_{r-k+1},\ ...,\ 
d_r \in p_r^{\theta_k(0,\delta_1, ..., \delta_{r-k-1})} c_kD_r ,
\end{equation}
where 
$$
D_j\ =\ \{p_j^i\ :\ |i| \leq I\},\ j=2,...,r,
$$
where $I$ is a suitably large integer.  A trivial upper bound for $I$ would be $r L^r \max_{i,j} |m_{i,j}|$.

Of course, with such a large collection of possible curves, some may fail to intersect
$(1/p_1,1] \times \cdots \times (1/p_r,1]$, and so will not contain any points of the form
(\ref{pointsn}) at all.  That's fine.  All we were interested in was a set of disjoint curves that
{\it do} cover all those points, and that can be suitably discretized later to prove certain theorems.

\subsubsection{An important property of the parameterized curves}
\label{importantproperty}

When all is said and done, the curves from the previous section that we generate have the
form
$$
(p_1^{t-1},\ e_2 p_2^{q_2 t},\ e_3 p_3^{q_3t},\ ...,\ e_r p_r^{q_r t}),
$$
where the $q_j$'s are non-zero rational numbers.  An important property here is the fact that
$p_1, p_2^{q_2}, ..., p_r^{q_r}$ are all distinct, which fulfills a goal mentioned at the
beginning of section \ref{paramsection}.  This property holds since if two of them were equal, 
we would have, for example,
$$
q_i \log p_i\ =\ q_j \log p_j,
$$
yet we know that the $\log p_i$'s are linearly independent over ${\mathbb Q}$.

\subsection{Discretized curves} \label{discretizedcurvessection}

Now we produce discretized versions of the curves produced in section \ref{paramsection}.  We begin
by defining ${\cal C}$ to be the set of all curves produced at the end of subsection \ref{paramsection}
with the property that the curve has non-empty intersection with the set 
$$
\Gamma\ :=\ (1/p_1,1]\ \times\ (1/p_2,1]\ \times\ \cdots\ \times\ (1/p_r,1].
$$
Now, any curve in ${\cal C}$ may be parameterized by a vector
\begin{equation}\label{paramvec}
(\delta_1, ..., \delta_{r-k-1}, c_1, ..., c_k, \tau_2, ..., \tau_r),
\end{equation}
where 
\begin{equation}\label{ctau}
c_i = p_{r-k+i}^{s_i},\ {\rm where\ }s_i \in S,\ i=1,...,k;\ {\rm and\ }\tau_i \in D_i,\ i=2, ..., r.
\end{equation}
The $\delta_1,...,\delta_{r-k-1}$ can take on a continuum of values in $[0,1)$, while the 
values taken on by the $c_1,...,c_k, \tau_2, ..., \tau_r$ are finite in number.  

Given a prime $P$ satisfying 
$$
P\ >\ (\max_{j=1,...,r} p_j)^{2H},
$$
we define a family of sets ${\cal F}$ as follows:  
for each choice of numbers $0 \leq f_1,...,f_{r-k-1} \leq P-1$, and
choice of $c_1,...,c_k, \tau_2, ..., \tau_r$ as above, let 
$F(f_1,...,f_{r-k-1},c_1,...,c_k,\tau_2,...,\tau_r)$ denote the set of all points
$$
(x_1,x_2, ..., x_r)\ \in\ \{0,1...,P-1\}^r,
$$
such that there exists a curve in ${\cal C}$ with parameter vector (\ref{paramvec}),
incident to a point $(y_1,...,y_r) \in \Gamma$,  such that
\begin{equation}\label{deltaf}
(\delta_1,...,\delta_r)\ \in\ (f_1/P, ..., f_r/P) + [0,1/P]^r,
\end{equation}
and
\begin{equation}\label{yi}
(y_1,...,y_r)\ \in\ (x_1/P,...,x_r/P) + [0,1/P]^r.
\end{equation}
If this set $F(f_1,...,f_{r-k-1},c_1,...,c_k,\tau_2,...,\tau_r)$ is non-empty, then we
add it to the family ${\cal F}$; otherwise, we don't.

One can easily see that, since there are at most $P^{r-k-1}$ choices
for $f_1,...,f_{r-k-1}$, and since there are only a bounded number
possibilities for the $c_j$'s and $\tau_j$'s, 
$$
|{\cal F}|\ \ll\ P^{r-k-1}.
$$
Likewise, for each choice of the $f_j$'s, there is at least one
choice of the other parameters making $F(f_1,...,f_{r-k-1}, c_1,...,c_k, \tau_2, ..., \tau_r)$ non-empty; and so, we have that
\begin{equation}\label{calFbound}
P^{r-k-1}\ \ll\ |{\cal F}|\ \ll\ P^{r-k-1}.
\end{equation}

\subsection{Two propositions and the proof of Theorem
\ref{theorem2}}

We assume that $P$ is as in the previous section.  
The two propositions we will need to prove Theorem \ref{theorem2}
are:

\begin{proposition}\label{propa}  We have that for 
every point $x \in \{0,1,2,...,P-1\}^r$, 
$$
\#\{n \leq N\ :\ (\alpha_1(n),..., \alpha_r(n)) \in x/P + 
[0,1/P)^r\}\ \ll\ N |{\cal F}|^{-1} P^{-1}.
$$
\end{proposition}

And:

\begin{proposition}\label{propb}  
Suppose $0 < \varepsilon < 1/20r^2$ and let
\begin{equation}\label{curveform}
K(t)\ =\ (\zeta_1 \theta_1^{t},\ \zeta_2 \theta_2^{t},\ ...,\ \zeta_r \theta_r^{t}),
\end{equation}
where the $\theta_1,...,\theta_r > 0$ are distinct real numbers
where none are equal to $1$, and where $\zeta_1,...,\zeta_r$ are 
also real (constants).  

\noindent Suppose that 
$$
F\ \subseteq\ \{0,1,...,P-1\}^r
$$
is the set of all vectors such that $(x_1,...,x_r) \in F$ if and 
only if there exists $t \in [0,1)$ such that 
\begin{equation}\label{KxP}
K(t)\ \in\ {1\over P}(x_1, x_2, ..., x_r) + [0,1/P)^r. 
\end{equation}

Now, let 
$$
A_1,\ A_2,\ ...,\ A_r\ \subseteq\ {\mathbb F}_P,\ {\rm with\ }|A_1|, ..., |A_r|\ \geq\ P^{1-\varepsilon}.
$$

We claim that for all but at most $o(P)$
elements $(x_1,...,x_r) \in F$, we have that for every $(\beta_1,...,\beta_r) \in {\mathbb F}_P^r$ 
there exist 
$$
1 \leq n \leq P^{10r^2\varepsilon},\ {\rm and\ } (\delta_1,...,\delta_r) \in \{0,1,...,[P^{7r\varepsilon}]\}^r,
$$
such that 
\begin{equation}\label{conclusion1}
n\cdot (x_1,...,x_r)\ +\ (\beta_1,...,\beta_r)\ -\ (\delta_1,...,\delta_r)\ \in\ (A_1+A_1+A_1) \times \cdots \times
(A_r + A_r + A_r).
\end{equation}
\end{proposition}

\noindent {\bf Remark:}  This proposition is near to "best-possible" since if we replace the bound on the $A_i$ with $|A_i| > P^{1-r \varepsilon}/3$
and replace the bound $n \leq P^{10r^2 \varepsilon}$ with $n \leq 
P^{r(r-1) \varepsilon/2}$, then the conclusion would not hold:  fix a number $\theta = 1/2$ and let $\theta_j = \theta^j$ for 
$j=1,2,...,r$, and let $\zeta_i = (-1)^{i-1} {r-1 \choose i-1}$.  
Choose $A_i$ to be all the integers in $[0,P^{1-r\varepsilon}/3]$.
Throughout, in all cases we will use 
$(\beta_1,...,\beta_r) = (0,...,0)$.

Now, we claim that for every $1 \leq n \leq P^{r(r-1) \varepsilon/2}$
and for every $t \in [0,1]$, if $(x_1,...,x_r) \in F$ satisfies 
(\ref{KxP}), then $n(x_1,...,x_r)$ cannot be an element of 
$3A_1 \times 3A_2 \times \cdots \times 3A_r$.  Note that, here, the $\delta_i$'s don't help much, 
since we could choose the $A_i$'s a little smaller, and even take the $\delta_i$'s of size up to $P^{1-r\varepsilon}$ or so, and still not get (\ref{conclusion1}).

To see this, we begin by noting that if $(x_1,...,x_r) \in F$ is 
associated to a particular $t$ in the sense of (\ref{KxP}), then
$$
\sum_{k=1}^r x_k\ \in\ 
P \sum_{k=1}^r (-1)^{k-1} {r-1 \choose k-1} \theta^{kt} - [0, r] + 
P\cdot {\mathbb Z}\ =\ P \theta^t (1 - \theta^t)^{r-1} - [0,r]
+ P \cdot {\mathbb Z}.
$$
Now, for $t \in [0,P^{-r\varepsilon}]$ we have that
\begin{eqnarray*}
P\theta^t (1-\theta^t)^{r-1}\ &\leq&\ P (1 - \exp(-t \log (1/\theta))^{r-1}
\ \leq\ P(1 - (1 - t \log(1/\theta)))^{r-1}\\
&<&\ P t^{r-1}\ <\ P^{1-r(r-1)\varepsilon}.
\end{eqnarray*}

So, for any $n \leq P^{r(r-1) \varepsilon/2}$, computing a dot-product
of $n(x_1,...,x_r)$ with $(1,...,1)$, we would have
$$
n(x_1+\cdots + x_r)
\ \in\ P^{1-r(r-1)\varepsilon/2} - [0,r] + P\cdot {\mathbb Z},
$$
which doesn't include any of the elements of $y_1+\cdots + y_r$, 
for $y_i \in 3A_i \subseteq [0, P^{1-r\varepsilon}]$.

\subsection{Completion of the proof of Theorem \ref{theorem2}}

We will later apply Proposition \ref{propb} using, for $j=1,...,r$,
\begin{eqnarray*}
A_j\ &:=&\ \left \{ d_1 \lceil P/p_j\rceil + d_2 \lceil 
P/p_j^2\rceil + \cdots + d_H \lceil P/p_j^H \rceil\ :\ 
0\leq d_1,...,d_H < p_j/10\right \}\\ 
&& \hskip0.5in +\ \{0,1,...,\lceil P/p_j^H\rceil -1\}.
\end{eqnarray*}

We note that 
$$
|A_j|\ \geq\ (p_j/10)^H (P/p_j^H)\ =\ P/10^H.
$$
This follows from the fact that all the expressions  
$$
d_1 \lceil P/p_j\rceil + \cdots + d_H \lceil P/p_j^H\rceil + x,
\ {\rm where\ } 0 \leq d_1,...,d_H < p_j/10,\ x \in \{0,1,...,
\lceil P/p_j^H\rceil-1\}
$$
are unique mod $P$
(it is easy to see this, using a similar proof as the
one showing base-$p_j$ representations are unique).  

We will have that, working in ${\mathbb F}_P$,
\begin{eqnarray}\label{3Asubset}
A_j+A_j+A_j &=& \{e_1 \lceil P/p_j\rceil + \cdots + 
e_H \lceil P/p_j^H\rceil\ :\ 0 \leq e_1,...,e_H \leq 3 \lfloor p_j/10\rfloor\}\nonumber \\
&&\hskip0.5in +\ \{0,1,...,3\cdot \lceil P/p_j^H\rceil -3\}
\end{eqnarray}
Thinking of this set as a subset of $\{0,1,2,...,P-1\}$, if 
we divide its elements by a factor $P$, then we get a set of numbers
contained in the set
$$
\left \{ {e_1 \over p_j} + 
{e_2 \over p_j^2} + \cdots + {e_H \over p_j^H}\ :\ 
0 \leq e_1,...,e_H < {3p_j \over 10}\right \} + \left [ 0,\ {4 \over p_j^H} \right ],
$$

We now let 
\begin{equation}\label{letF}
F\ =\ F(f_1, ..., f_{r-k-1}, c_1, ..., c_k, \tau_2, ..., \tau_r)\ \in
 {\cal F}
\end{equation}
be one of the sets in ${\cal F}$; $F$ is thus a discretized
version of a curve of general shape (\ref{curveform}),
where the $\zeta_j$'s depend on the choice of 
$f_1,..,f_{r-k-1}, c_1,...,c_k, \tau_2, ..., \tau_r$. 

We apply Proposition \ref{propb} to the curve $F$ using $\varepsilon$ defined by
$$
P^{1-\varepsilon}\ =\ \min_{j=1,...,r} |A_j|\ \geq\ {P \over 10^{H}};
$$
that is,
$$
P^\varepsilon\ \leq\ 10^H.
$$
Dividing (\ref{conclusion1}) through by a factor $P$ (interpreting coordinates now as integers instead of elements of ${\mathbb F}_P$), 
we get that for all but $o(P)$ of the $(x_1,...,x_r) \in F$, the
following holds:  let $y_1,...,y_r$ satisfy
$$
{x_j \over P}\ \leq\ y_j\ <\ {x_j + 1 \over P},
$$
and let $\beta_1,...,\beta_r \in [0,1)$ be arbitrary, and then
let $\beta'_1,...,\beta'_r \in \{0,1,...,P-1\}$ be the unique integers such that
for every $j=1,...,r$,
$$
{\beta'_j \over P}\ \leq\ \beta_j\ <\ {\beta'_j + 1 \over P}.
$$
Then, there exist
$$
1\ \leq\ s\ \leq\ P^{10r^2\varepsilon}\ <\ 10^{10r^2H}
$$
and 
$$
(\delta_1,..., \delta_r)\ \in\ \{0,1,...,
[P^{7r\varepsilon}]\}^r,
$$
such that for $j=1,...,r$,
\begin{eqnarray*}
s y_j + \beta_j\ \in\ {s x_j\over P} + {\beta'_j \over P} + 
\left [ 0,\ {1 + s \over P}\right ]\ &\subseteq&\ 
{\delta_j \over P} + {1\over P} \cdot (A_j+A_j+A_j) + \left [ 0,\ {1 + s \over P}\right ] + {\mathbb Z}\\
&\subseteq&\ U_j(H) + {\mathbb Z},
\end{eqnarray*}
where, recall, $U_j(H)$ is defined in (\ref{Uidef}).  
Note that in deducing this last containment 
we have used the fact that $0 < \varepsilon < 1/20r^2$, which
implies that 
$$
s/P,\ \delta_j/P\ <\ P^{10r^2 \varepsilon - 1}\ <\ P^{-1/2},
$$
which is much smaller than $1/p_j^H$, the width of the interval
in the definition of $U_j(H)$.  
Taking fractional parts of both sides, we get that, for all
$j=1,...,r$,
\begin{equation}\label{nyj}
\{sy_j + \beta_j\}\ \in\ U_j(H).
\end{equation}

We will use the notation
$$
F\ =\ F^\flat \sqcup F^\sharp,
$$
where $F^\sharp$ denotes the set of $x \in F$ such that
for every $\vec \beta$ there there exists $1 \leq s \leq P^{10 r^2 \varepsilon}$ where 
for every $y \in x/P + [0, 1/P)^r$ we have that 
(\ref{nyj}) holds; and we let $F^\flat$ denotes the
rest of $F$.  Note that from what we just proved, 
$|F^\flat| = o(P)$, and so $|F^\sharp| = |F| - o(P)$.
\bigskip

We will say that an integer $n \leq N$ is {\it good} if 
$$
\exists\ s \leq 10^{10r^2 H}\ \forall\ j =1,...,r,\ 
\{ s \alpha_j(n) + \beta_j(n)\} \in U_j(H)
$$
and, otherwise, we will say that it is {\it bad}.  We have that
the number of $n \leq N$ that are bad is at most 
$$
\sum_{F \in {\cal F}} \sum_{x \in F^\flat} 
\#\{n \leq N\ :\ (\alpha_1(n),...,\alpha_r(n)) 
\in x/P + [0,1/P)^r\}
$$
Applying Proposition \ref{propa} and (\ref{calFbound}) 
we get that this count is 
$$
\ll\ \sum_{F \in {\cal F}} N\cdot |F^\flat|\cdot |{\cal F}|^{-1} P^{-1}\ =\ o(P)\cdot N P^{-1}\ =\ o(N).
$$
This completes the proof of Theorem \ref{theorem2}.

\subsection{Proof of the Proposition \ref{propa}}

Fix a point $x=(x_1,...,x_r) \in {\mathbb F}_P^r$.  We will
only focus on counting the $n\leq N$ such that
$(\alpha_1(n), ..., \alpha_{r-k}(n))$ belongs to
$(x_1,...,x_{r-k})/P + [0,1/P)^{r-k}$.  This is legal, since
the proposition only claims an upper bound.

Now, since $\alpha_j(n) = p_j^{ \{n \log 2 / \log p_j\}-1}$, in
order for this to belong to $x_j/P + [0,1/P)$, we need that 
$\{n \log 2 / \log p_j\}$ belongs to a certain set $I_j + {\mathbb Z}$, where $I_j$ is an interval of width
at most $1/P \log p_j$.  Thus, our goal
is to count the number of $n \leq N$ such that
$$
\left ( \left \{ {n \log 2 \over \log p_1} \right \}-1,\ 
\left \{ {n \log 2 \over \log p_2}\right \}-1,\ ...,\ 
\left \{ {n \log 2 \over \log p_{r-k}}\right \}-1\right )\ 
\in\ I_1 \times I_2 \times \cdots \times I_{r-k} + {\mathbb Z}^{r-k}.
$$

Now, from (\ref{independentlogs}) and Theorem \ref{weyl} we
have that the number of such $n \leq N$ is, asymptotically,
$$
N (|I_1|\cdots |I_{r-k}| + o(1))\ \ll\ N P^{-r+k}\ \ll\ 
N |{\cal F}|^{-1} P^{-1}, 
$$
where the last expression follows from (\ref{calFbound}).
Note that the implied constants for the $\ll$'s depend on the 
$p_j$'s.  

This completes the proof since the upper bound on the set of
$n\leq N$ has the form claimed by the proposition.

\subsection{Proof of Proposition \ref{propb}}

Before embarking on the proof of this theorem, it's worth remarking
that although a more satisfying conclusion of the proposition would
omit the use of $(\delta_1,...,\delta_r)$, so that the conclusion
is something like
$$
n\vec x + \vec \beta\ \in\ 3A_1 \times 3A_2 \times \cdots \times 3A_r,
$$
we actually need to use this $(\delta_1,...,\delta_r)$ translate
in order that $q \geq 2$ in (\ref{qdef}) below.  Perhaps a more 
involved proof can get around the need for this translate, but given
that we apply a discretization process, passing from $K(t)$ to a 
set of points $F$, which can destroy some of the delicate arithmetic
properties of the curve $K(t)$, some care would be needed.

Also, we will assume that 
\begin{equation} \label{Flower}
|F|\ >\ P (\log P)^{-1},
\end{equation}
since otherwise the conclusion of the proposition is trivial.

\subsubsection{Basic setup}

For this proof we will use discrete Fourier methods.  Given a 
function $f : {\mathbb F}_P^r \to {\mathbb C}$, and a vector
$(a_1,...,a_r) \in \{0,1,2,...,P-1\}^r$,
we define the Fourier transform
$$
\widehat{f}(a_1,...,a_r)\ :=\ \sum_{(n_1,...,n_r) \in \{0,1,...,P-1\}^r} f(n_1,...,n_r) e^{2\pi i (a_1,...,a_r)\cdot (n_1,...,n_r)/P}.
$$

A consequence of Parseval is that
$$
\sum_{0 \leq s_1,...,s_r \leq P-1} |\widehat 1_{A_1 \times A_2 \times \cdots \times A_r}(s_1,...,s_r)|^2
\ =\ P^r |A_1|\cdots |A_r|.
$$
Thus, if $Q$ is the set of all places $(s_1,...,s_r)$ where
$$
|\widehat 1_{A_1 \times A_2 \times \cdots \times A_r}(s_1,...,s_r)|\ \geq\ P^{r(1 - 3\varepsilon)},
$$
then
$$
|Q|\ \leq\ P^{-2r(1-3\varepsilon)} P^r |A_1|\cdots |A_r|\ \leq\ P^{6r\varepsilon}  
$$

Let $Q' \subseteq Q$ be all those places $(s_1,...,s_r) \in Q$
satisfying
\begin{equation}\label{sbounds}
|s_i|\ \leq\ P^{1-6r\varepsilon},\ i=1,2,...,r.
\end{equation}
Let $N = |Q'|$, and note that
\begin{equation} \label{Nbounds}
N\ \leq\ |Q|\ \leq\ P^{6r\varepsilon}.
\end{equation}

Now we let $E$ denote the set of all $(x_1,...,x_r)\in F$, such that  there exists
$(s_1,...,s_r) \in Q'$, $(s_1,...,s_r) \neq (0,...,0)$, such that
\begin{equation}\label{Eineq}
\left \|{(x_1,...,x_r) \cdot (s_1,...,s_r) \over P} \right \|\ =\ 
\left \|{ x_1 s_1 + \cdots + x_r s_r\over P}\right \|\ <\ 
{1 \over P^{8r^2\varepsilon}}.
\end{equation}

\subsubsection{Proposition follows if we can show $|E|=o(|F|)$}

We will show that $|E| = o(|F|)$.  If this holds, then let us see how it implies the conclusion of the
Proposition:  let $L = [\log P]$, 
$$
U\ :=\ \{0,1,2,...,[P^{7r\varepsilon}/L]\}^r,
$$
and define
$g(\vec \delta) = g(\delta_1,...,\delta_r)$ to be the following $L$-fold convolution
$$
g(\vec \delta)\ :=\ 1_U * 1_U * \cdots *1_U(\delta_1,...,\delta_r).
$$
Now, let 
\begin{equation}\label{xset}
(x_1,...,x_r)\ \in\ F \setminus E
\end{equation}
be any of the $|F|-o(|F|)$ vectors such that
(\ref{Eineq}) fails to hold, for every $(s_1,...,s_r) \in Q'$.
Let 
$$
M\ :=\ [P^{10 r^2 \varepsilon}],
$$
and let $f$ be the indicator function for the set
$$
\{ (-nx_1,-nx_2, ..., -nx_r)\ :\ 1 \leq n \leq M\}.
$$

Then, we have that if
\begin{equation}\label{target1}
1_{A_1 \times \cdots \times A_r} * 1_{A_1 \times \cdots \times A_r} * 1_{A_1 \times \cdots \times A_r} * g * f(\vec \beta)\ >\ 0,
\end{equation}
then there exists $1 \leq n \leq M$ and $(\delta_1,...,\delta_r)$,
so that (\ref{conclusion1}) holds.  

Expressing the left-hand-side of (\ref{target1}) in terms of 
Fourier transforms, one sees that it equals:
\begin{eqnarray}\label{target2}
&& P^{-r} \sum_{(s_1,...,s_r) \in {\mathbb F}_P^r}
e^{-2\pi i \vec s \cdot \vec \beta/P} \widehat 1_{A_1\times \cdots \times A_r}(s_1,...,s_r)^3
\widehat g(s_1,...,s_r) \widehat f(s_1,...,s_r) \nonumber \\
&&\hskip0.5in =\ P^{-r} \sum_{\vec{s} \in {\mathbb F}_P^r}
e^{-2\pi i \vec s \cdot \vec \beta/P}\widehat 1_{A_1\times \cdots \times A_r}(\vec{s})^3 
\widehat 1_U(\vec{s})^L \widehat f(\vec{s}).
\end{eqnarray}

We split the terms in the second sum into the term with $(s_1,...,s_r) = (0,...,0)$,
the terms $(s_1,...,s_r) \in Q$, and then the remaining terms.

The contribution of the term $(s_1,...,s_r) = (0,...,0)$ is
\begin{equation}\label{zerocount}
P^{-r} M |U|^L |A_1|^3 \cdots |A_r|^3.
\end{equation}

Now suppose $(s_1,...,s_r) \in Q \setminus Q'$.  Then, for
some $i=1,...,r$ we have that $P^{1-7r\varepsilon} < |s_i| < P/2$.
Thus,
$$
|\widehat g(\vec s)|\ \ll\ \prod_{i=1}^r 
\min(|U|^{L/r}, \| s_i/P\|^{-L})
\ <\ |U|^{L(r-1)/r} P^{6r\varepsilon L}\ \leq\ 
|U|^L P^{-rL\varepsilon}.
$$
It follows, then, that the contribution of all such 
$(s_1,...,s_r) \in Q \setminus Q'$ to the right-hand-side of
(\ref{target2}) is bounded from above by 
$$
P^{-r} N |A_1|^3\cdots |A_r|^3 |U|^L P^{-rL\varepsilon} M,
$$
which is much smaller than (\ref{zerocount}), on account of the 
$P^{-rL\varepsilon}$ factor, even when using the crude 
upper bound $N \leq P^{6r\varepsilon}$.  

Next, we consider the contribution of all terms with
$(s_1,...,s_r) \in Q'$.  Then, since $(x_1,...,x_r)$
satisfies (\ref{xset}), and in particular that it is not $E$, 
we have that 
\begin{eqnarray*}
|\widehat f(s_1,...,s_r)|\ &=&\ \left | \sum_{1 \leq n \leq M} 
e^{2\pi i n (x_1,...,x_r) \cdot (s_1,...,s_r)/P} \right |\\
&\ll&\ {1 \over 
\| (x_1,...,x_r)\cdot (s_1,...,s_r)/P\|}\\
&\leq&\ P^{8r^2\varepsilon}.
\end{eqnarray*}
So, the contribution of the terms in (\ref{target2}) with $(s_1,...,s_r) \in Q'$, 
$(s_1,...,s_r) \neq (0,...,0)$, is, by Parseval,
\begin{eqnarray*}
&\ll&\ P^{-r} P^{8r^2\varepsilon} 
|U|^L \sum_{0\leq s_1,...,s_r \leq P-1} 
|\widehat 1_{A_1\times \cdots \times A_r}(s_1,...,s_r)|^3\\
&\leq&\ P^{-r + 8r^2\varepsilon}|U|^L |A_1| \cdots |A_r| 
\sum_{0\leq s_1,...,s_r \leq P-1} |\widehat 1_{A_1\times \cdots \times A_r}(s_1,...,s_r)|^2 \\
&\leq&\ P^{8r^2\varepsilon} |U|^L |A_1|^2 \cdots |A_r|^2 \\
&\ll&\ P^{-r-r\varepsilon} M |U|^L |A_1|^3 \cdots |A_r|^3,
\end{eqnarray*}
which is smaller than the contribution of the term with $(s_1,...,s_r) = (0,...,0)$ given in 
(\ref{zerocount}).

Finally, we consider the contribution of the remaining terms.  For these terms we have
$$
|\widehat 1_{A_1 \times \cdots \times A_r}(s_1,...,s_r)|\ <\ P^{r(1-3\varepsilon)}\ \leq\ 
|A_1|\cdots |A_r| P^{-2r\varepsilon}.
$$
Using this in those terms on the right-hand-side of (\ref{target2}), we find that,
using Parseval again, they contribute at most
\begin{eqnarray*}
&& P^{-r-2r\varepsilon} M |U|^L |A_1|\cdots |A_r|\sum_{0\leq s_1,...,s_r \leq P-1} |\widehat 1_{A_1\times \cdots \times A_r}(s_1,...,s_r)|^2\\
&&\ \ \ \ \leq\ P^{-2r\varepsilon} M |U|^L
|A_1|^2 \cdots |A_r|^2\ \leq\ 
P^{-r - r\varepsilon} M |U|^L |A_1|^3\cdots |A_r|^3,
\end{eqnarray*}
which is also appreciably smaller than the contribution of the term with
$(s_1,...,s_r) = (0,...,0)$, as in (\ref{zerocount}).

Thus, there exists $1 \leq n \leq M$ and 
$0 \leq \delta_1,...,\delta_r \leq P^{7r\varepsilon}$ 
so that 
$$
n(x_1,...,x_r) +\vec \beta - \vec \delta\ \in\ (3A_1) \times (3A_2) \times \cdots \times (3A_r).
$$
And since this holds for $(1-o(1))|F|$ vectors $(x_1,...,x_r) \in F$, the proposition is proved.

\subsubsection{Proving $|E|=o(|F|)$}

We begin by noting that we may assume that $Q'$ contains at least
one non-zero vector, since otherwise in the previous subsection
we never need to make use of bounds on $|\widehat f(s_1,...,s_r)|$,
nor reference to $(x_1,...,x_r)$ -- we obtain the same bounds
independent of choice of $(x_1,...,x_r)$, which would imply that
$E$ is empty.

We note, by the pigeonhole principle, that 
there exist  $(s_1,...,s_r) \in Q'$, 
such that (\ref{Eineq}) holds for at least 
$|E|/N$ vectors $(x_1,...,x_r) \in E$.  Call this new set of 
vectors $E' \subseteq E$; so, we have
$$
|E'|\ \geq\ |E| / N.
$$

Let us suppose, for proof by contradiction, that
\begin{equation}\label{contradictme}
|E|/N\ >\ |F| P^{-7r\varepsilon} (\log P)^3.
\end{equation}
Note that if we establish a contradiction, then we would be forced to conclude that
$$
|E|\ \leq\ N |F| P^{-7r\varepsilon} (\log P)^3\ \leq\ 
|F| P^{-r\varepsilon+o(1)},
$$
which would imply $|E| = o(|F|)$, and which is just what we wanted to show.  

It may seem like we are throwing away a lot by applying the Pigeonhole Principle in this way, but at the end of this subsection we will give a plausibility argument for why the approach is, in fact, optimal in some sense.
\bigskip

For each $\vec x = (x_1,...,x_r) \in E'$, let 
$t = t(\vec x)$ be any value of $t$, so that
\begin{equation}\label{yvec}
K(t)\ \in\ (x_1/P, ..., x_r/P) + [0,1/P)^r.
\end{equation}
Also, for any vector $\vec v \in [0,1)^r$, let $\pi(v)$ denote the
unique $\vec x \in \{0,...,P-1\}^r$, so that
$$
\vec v\ \in\ {1\over P} \vec x + \left [0,\ {1\over P}\right ]^r.
$$
Note that we can insist on the set of $t : E' \to {\mathbb R}$ be
injective, if need be.

Now, if we consider the set of all points in a cube
\begin{equation}\label{w_cube}
\vec w + \left [0,{1 \over P}\right ]^r,
\end{equation}
where $\vec w$ is some arbitrary $r$-dimensional vector, the
function $\pi$ will map that set to a set of size at most $2^r$.
Thus, if we let 
$$
\Upsilon\ :=\ \max_{i=1,...,r} |\zeta_i|,
$$
and let
$$
T\ :=\ \{t(\vec x)\ :\ \vec x \in E'\},
$$
then we claim that any interval of
width $(\Upsilon P)^{-1}$ can have at most $2^r \log P$ elements of 
$T$.  

The reason this holds is that if we restrict $t$ to an interval $I$ of width at most $(\Upsilon P \log P)^{-1}$,
then the coordinates of $K(t)$ will vary by $o(1/P)$; and so,
the set $\{ K(t)\ :\ t\in I\}$ will be contained in one of
the cubes (\ref{w_cube}), which can correspond to at most 
$2^r$ vectors $\vec x \in \{0,1,...,P-1\}^r$.

By picking at most one element of $T$ in each interval of width 
$(\Upsilon P)^{-1}$, we can pass to a subset
$$
T'\ \subseteq\ T,\ {\rm where\ } |T'| > 2^{-r} |T| (\log P)^{-1} = 
2^{-r} |E'| (\log P)^{-1} > |F| P^{-7r\varepsilon}(\log P)^2,
$$
such that every pair of elements of $T'$ is at least $(\Upsilon P)^{-1}$ apart.

Now we index the elements or $T'$ as follows:
$$
T'\ :=\ \{t_1,t_2, ..., t_n\}, 
$$
where 
$$
t_1\ <\ t_2\ <\ \cdots\ <\ t_n.
$$

Then, we extract disjoint subsets $T_1,...,T_{2^r} \subseteq T'$ as
follows:  we let
$$
T_i\ :=\ \{t_j\ :\ (2i-2)n/2^{r+1} < j < (2i-1)n/2^{r+1}\},
$$
which satisfies
\begin{equation}\label{Tibound}
|T_i|\ \gg\ n/2^r\ \gg\ |T'|\ >\ |F| P^{-7r\varepsilon} (\log P)^2.
\end{equation}
Let 
$$
d(T_i,T_j)\ :=\ \min_{t \in T_i, u \in T_j} |t-u|.
$$
Since the elements of $T'$ are spaced at least $(\Upsilon P)^{-1}$ 
apart, using (\ref{Flower}) we must have that  
\begin{equation}\label{separation}
\min_{1 \leq i < j \leq 2^r} d(T_i,T_j)\ \geq\ n/2^{r+1}\Upsilon P\ \gg\ n / |F|\ >\ P^{-7r\varepsilon} (\log P)^2,
\end{equation}
which follows from the fact that if $i$ is the coordinate where
$\Upsilon = |\zeta_i|$, then for each change in $t$ by $\asymp (\Upsilon P)^{-1}$, $\zeta_i \theta_i^t$ changes by at least $1/P$.  

Define, also, the associated intervals
$$
I_i\ :=\ [t_{\lceil (2i-2)n/2^{r+1}\rceil},\ t_{\lfloor (2i-1)n/2^{r+1}\rfloor}]. 
$$
Note that if $T_i \subset I_i$.

We now define $u_1,...,u_{2^r}$ as follows:  we let $u_i$ be any 
element in the interval $I_i$ such that $|h'(u)|$ is minimal, where
$$
h(t)\ :=\ (s_1,...,s_r) \cdot K(t)\ =\ s_1 \zeta_1 \theta_1^{t} + \cdots + s_r \zeta_r \theta_r^{t}.
$$
Note that
\begin{equation} \label{h'te}
h'(t)\ :=\ s_1 \zeta_1 \theta_1^t \log \theta_1 + s_2 \zeta_2 \theta_2^t \log \theta_2  + \cdots + 
s_r \zeta_r \theta_r^t \log \theta_r.
\end{equation}

We now need the following lemma, which makes use of an idea from
\cite[page 99, book 2, example 1]{gantmacher}:
\bigskip

\begin{lemma}  Suppose
$$
H(t) = c_1 x_1^t + \cdots + c_r x_r^t,\ {\rm where\ }x_1,...,x_r > 0
\ {\rm distinct}.
$$
Then, for any sequence 
$$
0\ <\ v_1\ <\ v_2\ <\ \cdots\ <\ v_{2^r}\ <\ 1,
$$
there exists $j=1,...,2^r$ such that 
$$
|H(v_j)|\ \geq\ \Delta^{r-1} c(x_1,...,x_r) \max_{i=1,...,r} |c_i|,
$$
where $\Delta$ is the minimum of the differences 
$|v_{k+1}-v_k|$, $k=1,...,2^r-1$, and $c(x_1,...,x_r)$ is
some constant that depends only on $x_1,...,x_r$.
\end{lemma}

\noindent {\bf Remark:}  It may be possible to obtain the conclusion using fewer than $2^r$ vector $v_i$.  However, the factor $\Delta^{r-1}$ in the bound for $|H(v_j)|$ in the lemma is essentially best-possible for the following reason:  consider the case where
$$
H(t)\ =\ (e^t - 1)^{r-1}\ =\ \sum_{j=0}^{r-1} (-1)^j {r-1 \choose j} (e^j)^t. 
$$
We can let $x_j = e^{j-1}$, $j=1,...,r$, and then let 
$c_j = (-1)^{j-1} {r-1 \choose j-1}$.  Then, if we choose 
$v_j = j \Delta$ we will see that, for $\Delta$ taken as small 
as desired,
$$
|H(v_j)|\ =\ (e^{v_j}-1)^{r-1}\ =\ (v_j + O(v_j^2))^{r-1}\ \ll\ 
j^{r-1} \Delta^{r-1}.
$$
\bigskip

\noindent {\bf Proof.}  This can be easily proved by 
induction.  First note that the claim holds when $r=1$.
Assume, for proof by induction, the claim holds for 
$r=\ell$.  Consider now the case $r=\ell+1$:  we first
$i_0=1,...,r$ be such that $|c_{i_0}|$ is maximal.  We then let 
$j_0=1,...,r$ be any index so that $j_0 \neq i_0$; in fact, just set it
to $2$ if $i_0=1$ and set it to $1$ is $i_0 \geq 2$.  

Then we define
$$
K(t)\ =\ x_{j_0}^{-t} H(t)\ =\ c_{j_0} + \sum_{1 \leq k \leq \ell+1 \atop 
k \neq j_0} c_k (x_k/x_{j_0})^t.
$$
By the Mean Value Theorem we have that there exists a sequence
$$
v'_j \in [v_{2j-1}, v_{2j}],\ j=1,...,2^{r-1} = 2^\ell,
$$
such that 
\begin{equation}\label{inductwith}
|K(v_{2j}) - K(v_{2j-1})|\ =\ |v_{2j} - v_{2j-1}|\cdot 
|K'(v'_j)|\ \geq\ \Delta |K'(v'_j)|.
\end{equation}
Now, 
$$
K'(v'_j)\ =\ \sum_{1 \leq k\leq \ell+1 \atop k \neq j_0} c_k \log(x_k/x_{j_0})
(x_k/x_{j_0})^t,
$$
which involves $\ell$ terms; so, we can apply the induction
hypothesis to it.  Note first that
$$
\min_{j=1,...,2^\ell-1} |v'_{j+1}-v'_j|\ \geq\ 
\min_{j=1,...,2^\ell} |v_{2j}-v_{2j-1}|\ \geq\ \Delta.
$$
So, now, applying the induction hypothesis to $K'$, we have
that for some $j=1,...,2^\ell$,
$$
|K'(v'_j)|\ \geq\ \Delta^{\ell-1} c(x_1/x_{j_0}, ..., x_{j_0-1}/x_{j_0},
x_{j_0+1}/x_{j_0}, ..., x_{\ell+1}/x_{j_0}) \max_{1\leq i \leq \ell+1
\atop i\neq j_0} |\log(x_i/x_{j_0})c_i|.
$$
Note that the $c(\cdots)$ here is a function of $\ell$ terms, not 
$\ell+1$ like we aim to show for $H(t)$ using $r=\ell+1$.  
Now, we can replace the maximum over $i$ of $|\log(x_i/x_{j_0})c_i|$
by something smaller; and it's clear we can take this smaller thing to 
be
$$
\max_{1\leq i \leq \ell+1} |c_i| 
\min_{1\leq i \leq \ell+1\atop i\neq j_0} |\log(x_i/x_{j_0})|\ =\ 
|c_{i_0}| \min_{1\leq i \leq \ell+1\atop i\neq j_0} |\log(x_i/x_{j_0})|.
$$
So, we have a bound of the shape
$$
|K'(v_j')|\ \geq\ \Delta^{\ell-1} c_2(x_1,...,x_{\ell+1}) |c_{i_0}|,
$$
where
$$
c_2(x_1,...,x_{\ell+1})\ =\  c(x_1/x_{j_0}, ..., x_{j_0-1}/x_{j_0},
x_{j_0+1}/x_{j_0}, ..., x_{\ell+1}/x_{j_0}) \min_{1\leq i \leq \ell+1\atop i\neq j_0} |\log(x_i/x_{j_0})|.
$$

And so, from (\ref{inductwith}) we conclude
$$
|K(v_{2j}) - K(v_{2j-1})|\ \geq\ \Delta^\ell c_2(x_1,...,x_{\ell+1})
|c_{i_0}|.
$$
Thus, either $|K(v_{2j})|$ or $|K(v_{2j-1})|$ is at least as big as
the right-hand-side.  It follows that  
\begin{eqnarray*}
\max (|H(v_{2j})|,\ |H(v_{2j-1})|)\ &\geq&\ 
\min(1, x_{j_0}) \max (|K(v_{2j})|,\ |K(v_{2j-1})|) \\
&\geq&\ \min(1,x_{j_0}) \Delta^\ell c_2(x_1,...,x_{\ell+1}) |c_{i_0}| \\
&\geq&\ \Delta^\ell  c_3(x_1,...,x_{\ell+1}) |c_{i_0}|,
\end{eqnarray*}
where
$$
c_3(x_1,...,x_{\ell+1})\ =\ c_2(x_1,...,x_{\ell+1}) \min(1, x_{j_0}).
$$
This then proves the induction step, since $c_3(\cdots)$ is a function
purely of $x_1,...,x_{\ell+1}$ as claimed by the lemma.
\hfill $\blacksquare$
\bigskip

Applying the Lemma to (\ref{h'te}), and also using (\ref{separation}),
and letting 
$$
\Delta\ :=\ \min_{1\leq i < j \leq 2^r} |u_i-u_j|\ \gg\ P^{-7r\varepsilon} (\log P)^2,
$$
we find that there exists
$i=1,...,2^r$ such that
\begin{eqnarray}\label{h'ui}
|h'(u_i)|\ &\geq&\ \Delta^{r-1} c(\theta_1,...,\theta_r) 
\max_j |s_j \zeta_j \log \theta_j|\nonumber \\
&\geq&\ P^{-7r(r-1)\varepsilon} 
c(\theta_1,...,\theta_r) \max_j |s_j \zeta_j \log \theta_j|.
\end{eqnarray}
And, thus, we will have the same bound for 
$|h'(t)|$ for every $t\in I_i$, by the way we chose $u_i$ minimally.  

Note, in particular, this means that $h'(t) \neq 0$ for $t \in I_i$;
so, $h(t)$ is either increasing on all of $I_i$ or decreasing on all of
$I_i$.  Without loss, let us assume that it is increasing on $I_i$.  

Now, we have from the triangle inequality and (\ref{Flower})
that for every $t \in I_i$,
$$
|h(t)|\ \leq\ B\ :=\ r \cdot c'\cdot \max_{j=1,...,r}
|s_j \zeta_j|\ \leq\ r \cdot c'\cdot P^{1-7r\varepsilon} \Upsilon
\ \ll\ |F| P^{-7r\varepsilon} \log P,
$$
where $c' = \max(1,\theta_1,...,\theta_r)$.  Thus, for each $t \in I_i$,
the nearest integer to $h(t)$ lies in the set 
$\{-B,-B+1,...,0,...,B\}$.

Now, by the pigeonhole principle we have that there exist an 
integer $z$ such that there are at least
\begin{equation}\label{qdef}
q\ :=\ \left \lfloor 
{|T_i| \over 2B+1}\right \rfloor\ \gg\ \log P.
\end{equation}
values $t \in T_i$ where $z$ is the nearest integer to $h(t)$.  Let $T'_i$ denote this set of values $t \in T_i$.  
\bigskip

We now want to see that the above is impossible:  let those $q$ values
of $t\in T'_i$ be
$$
t'_1\ <\ t'_2\ <\ \cdots\ <\ t'_q.
$$
From (\ref{Eineq}) we know that for $t \in T'_i$,
$$
\| h(t)\|\ \ll\ P^{-8r^2 \varepsilon},
$$
and so it follows that
\begin{equation}\label{h'iz}
h(t'_j)\ =\ z + O(P^{-8r^2\varepsilon}),\ {\rm for\ }j=1,...,q.
\end{equation}

Now, since $h$ is increasing on $I_i$ we have 
$$
h(t'_1)\ <\ h(t'_2)\ <\ \cdots\ <\ h(t'_q).
$$
But we also have by the Mean Value Theorem, and the fact
that $t_j > t_{j-1} + 1/\Upsilon P$ and (\ref{Flower}),
\begin{eqnarray*}
h(t'_j) - h(t'_{j-1})\ \geq\ (t'_j - t'_{j-1}) h'(u_j)
\ &\geq&\ \Upsilon^{-1} P^{-1-7r(r-1)\varepsilon} c(\theta_1,...,\theta_r)
\max_\ell |s_\ell \zeta_\ell \log \theta_\ell| \\
&\gg&\ (|F| \log P)^{-1} P^{-7r(r-1)\varepsilon} 
\max_\ell |s_\ell \zeta_\ell|.
\end{eqnarray*}
So, by telescoping, and applying (\ref{Tibound}), we find that
\begin{eqnarray}\label{finalcontradiction}
h(t'_q) - h(t'_1)\ &\gg&\ q (|F| \log P)^{-1} P^{-7r(r-1)\varepsilon} \cdot 
\max_j |s_j \zeta_j| \nonumber \\
&\gg&\ {|T_i|\cdot \max_j |s_j \zeta_j| \over 
B |F| P^{7r(r-1)\varepsilon} \log P} \nonumber \\
&\geq&\ P^{-7r^2 \varepsilon - o(1)}. 
\end{eqnarray}
(Note that here we have secretly used the fact that all the 
$\theta_j \neq 1$, which ensures that the $\min_j |\log \theta_j|
\neq 0$.)

But this clearly contradicts (\ref{h'iz}).  We conclude, therefore, that
$|E| = o(|F|)$, as needed.
\bigskip

Lastly, as promised at the beginning of this subsection, we explain why passing to a single $(s_1,...,s_r) \in Q'$, and then attempting to derive a contradiction to (\ref{contradictme}) is in some sense an optimal
approach -- in the sense that it seems unlikely we could easily show
$|E|$ is much smaller, say
\begin{equation}\label{Ecant}
|E|\ =\ o(P^{1-2r\varepsilon}).
\end{equation}
Note that in the above argument there is a little slack, as the 
$7$ at the end of (\ref{finalcontradiction}) is smaller than the $8$ in 
(\ref{h'iz}); and if we tried to eliminate this slack, instead of
deducing $|E| \ll P^{1-r\varepsilon}$, we would have a bound closer
to $|E| \ll P^{1-2r\varepsilon}$ -- so, if we found a heuristic argument
for why (\ref{Ecant}) cannot hold, it would mean that the bounds in our approach to showing $|E| = o(P)$ are near to best-possible.

Let us suppose we have that
$$
K(t)\ =\ (2^t,2^{2t},2^{3t}, ..., 2^{rt}),
$$
so that $\theta_i=2^{it}$, $i=1,...,r$, and that $\zeta_i = 1$.
For $j=1,...,N-1$, let $\nu_j = j \Delta$, where we take
$\Delta = 1/N$.  Then consider the function
$$
f_j(t)\ :=\ 2^t(2^t - 2^{\nu_j})^{r-1}\ =\ \sum_{\ell=0}^{r-1} 
{r-1 \choose \ell} (-2^{\nu_j})^{r-1-\ell} (2^{\ell+1})^t.
$$
For $t \in \nu_j + [-\Delta P^{-2r\varepsilon},\Delta P^{-2r\varepsilon}]$ we have that
$$
|f_j(t)|\ \ll\ (\Delta P^{-2r\varepsilon})^{r-1}\ =\ N^{-(r-1)} P^{-2r(r-1)\varepsilon},
$$
and all we know about $N$ is the upper bound (\ref{Nbounds}), so in
the worst case (largest $N$) we would have
\begin{equation}\label{fjtbound}
|f_j(t)|\ \ll\ P^{-8r(r-1) \varepsilon}.
\end{equation}

Note that the union $\cup_{j=1}^N (\nu_j + [-\Delta P^{-2r\varepsilon}, \Delta P^{-2r\epsilon}]$ has measure $\asymp P^{-2r\epsilon}$; and so,
we would expect that the number of $\vec x \in F$ with 
$t(\vec x)$ in this union to be $\asymp P^{1-2r\varepsilon}$.

And, we may interpret $f_j(t)$ as a certain dot-product -- it's the
dot-product of $K(t)$ with the vector whose $i$th coordinate is
${r-1\choose \ell} (-2^{\nu_j})^{r-1-\ell}$.  For each $j=1,...,N-1$
let us name this vector $V_j$.  Thus,
$$
f_j(t)\ =\ K(t)\cdot V_j.
$$
Setting aside for the time being that these vectors $V_j$ don't
necessarily have integer coordinates, if we had that $Q'$ contained all
of them (recall $|Q'| = N$, so it's large enough to contain them),
then we would have something like that for every 
$(x_1,...,x_r)\in F$ corresponding to $t$ in the above union, there exists $j=1,...,N-1$ such that
$$
\left \| {(x_1, ..., x_r) \cdot V_j \over P} \right \|\ \ll\ {1 \over P^{8r(r-1)\varepsilon}},
$$
which is close to what we see in (\ref{Eineq}).  

\section{Acknowledgements}

We would like to thank Thomas Bloom for some helpful conversations about this work.

\end{document}